\voffset-1cm
\input amstex
\documentstyle{amsppt}  
\mag=\magstep1  
\NoRunningHeads  
\NoBlackBoxes
\define\dist{\operatorname{dist}}
\define\gr{\operatorname{gr}}
\define\id{\operatorname{id}}
\define\rel{\operatorname{rel}}
\topmatter 
\title Boundary entropy of a hyperbolic group\endtitle  
\author Andrzej Bi\'s and Pawe\l \/ Walczak\endauthor  
\address Wydzia\l\/ Matematyki, Uniwersytet \L\'odzki, Banacha 22, 90238 
\L\'od\'z,  Poland\endaddress  
\email andbis\@ imul.uni.lodz.pl and pawelwal\@ imul.uni.lodz.pl \endemail  
\date November 24, 1998\enddate
\thanks The second author was supported by KBN grant 1037/P03/96/10.\endthanks
\abstract We show that the entropy of a hyperbolic group acting on its ideal 
boundary is closely related to the exponential rate of its growth.\endabstract
\subjclass 20 F 32\endsubjclass  
\endtopmatter  
 
\document 
\subhead 0. Introduction\endsubhead 
Hyperbolic groups in the Gromov's sense [{\bf Gro}] play an important role
in geometric group theory (see [{\bf GrH}] and the references there). In
particular, any non-elementary hyperbolic group has exponential growth and
the compact boundary of positive finite Hausdorff dimension ([{\bf GhH}],
pp. 126 and 157). Also, the group itself acts on its boundary {\it via}
Lipschitz quasi-conformal maps (ibidem, p. 127). The dynamics of this 
action is of great interest. For instance, it has been shown [{\bf CP}]
that this action is finitely presented, i.e. it is semiconjugate 
to a subshift of finite type in such a way that the fibres of the conjugating
map are finite of bounded length and the equivalence relation determined
by this map (two points are related whenever their images are equal)
is  another subshift of finite type. Also, one 
can consider the topological entropy of this action in the sense of [{\bf
GLW}]. In this article we prove the following.

\proclaim{Theorem} The topological entropy (with respect to a finite 
symmetric generating set) of a hyperbolic group $G$ acting on its ideal
boundary lies between the  exponential rate of growth of $G$ relative
to suitable bounds depending on the geometry of the group and
the exponential rate of growth of $G$ (with respect to the same
generating set). 
\endproclaim

The precise description of the bounds mentioned in the Theorem can be found
in Section 4 which contains also the proof of the Theorem and some final
remarks. In the first section, we recall the notion of the topological entropy
of a group action. In Section 2, we define the exponential rate of growth
relative to given constants. In Section 3, we provide a short review on
hyperbolic groups and spaces.

\subhead 1. Entropy\endsubhead
Let $G$ be a finitely generated group of homeomorphisms of a compact
metric space $(X, d)$ and $S$ be a finite symmetric ($e\in S$, $S^{-1} =
S$) set generating $G$. Equip $G$ with the word metric $d_S$ induced by
$S$ and let $B(n)$, $n\in\Bbb N$, denote the ball in $G$ of radius $n$ 
and centre $e$. Two points $x$ and $y$ of $X$ are said to be
$(n, \epsilon )$-separated ($\epsilon > 0$, $n\in\Bbb N$) whenever
$$d(gx, gy)\ge\epsilon$$
for some $g\in B(n)$. Since $X$ is compact, the maximal number
$N(n,\epsilon)$ of pairwise $(n, \epsilon )$-separated points of $X$ is
finite. Also, there exist finite $(n, \epsilon )$-spanning subsets of $X$:
A subset $A$ of $X$ is $(n,\epsilon )$-spanning whenever for any $y\in X$
there exists $x\in A$ such that
$$d(gx, gy) < \epsilon$$
for all $g\in B(n)$. Let $N'(n, \epsilon )$ denote the minimal cardinality
of an $(n,\epsilon )$-spanning subset of $X$.

Similarly to the case of classical dynamical systems ([{\bf Wa}], p. 169),
the families $N(n, \epsilon )$ and $N'(n,\epsilon )$ of functions have the
same type of growth [{\bf Eg}], more precisely, they have the same rate of
exponential growth and the topological entropy $h(G, S)$ of $G$ (w.r.t.
$S$) can be defined by the formula
$$h(G, S) = \lim_{\epsilon\to 0}\limsup_{n\to\infty}\frac{1}{n}\log N(n,
\epsilon ) =  \lim_{\epsilon\to 0}\limsup_{n\to\infty}\frac{1}{n}\log N'(n,
\epsilon ).$$

If all the maps of $G$ are Lipschitz and $X$ has finite Hausdorff
dimension, then $h(G, S)$ is finite for any $S$ (compare [{\bf GLW}],
Prop.2.7). Also, if $h(G, S) = 0$ for some $S$, then $h(G, S') = 0$ for any 
other generating set $S'$. Therefore, one can distinguish between groups of
positive and vanishing entropy without referring to generating sets.

\subhead 2. Growth\endsubhead
Let us keep the notation of the previous section and recall that the 
exponential {\it rate of growth} of $G$ (with respect to $S$) is defined as
$$\gr (G, S) = \lim_{n\to\infty}\frac{1}{n}\log N(n) 
= \lim_{n\to\infty}\frac{1}{n}N_0(n),$$
where $N(n) = \# B(n)$ and $N_0 (n) = \# S(n)$ is the cardinality of the sphere
$S(n)$ of radius $n$ and centre $e$. Also, if $\epsilon > 0$ and 
$N_0(n; \epsilon )$ is the maximal cardinality of an $\epsilon$-separated
subset $A$ of $S(n)$, then
$$\gr (G,S) = \lim_{n\to\infty}\frac{1}{n}\log N_0(n; \epsilon ).$$
In fact, $N_0(n; \epsilon )\le N_0(n)$ and, if $A$ is such a subset
of $S(n)$, then $\cup_{x\in A}B(x,\epsilon )\supset  S(n)$ and, therefore,
$N_0(n; \epsilon )N(\epsilon )\ge N_0(n)$ for any $n\in\Bbb N$.

Moreover, if $m\in\Bbb N$ and $A_n$ is a maximal $\epsilon$-separated subset
of $S(mn)$, $n = 1, 2, \dots$, then for any $x\in A_n$ we can find a sequence
$(x_0, x_1, \dots x_n)$ of elements of the group $G$ for which $x_k\in A_k$,
$x_n = x$, $x_0 = e$ and  $d(x_k, x_{k+1})\le m + \epsilon$ for all $k$. 
To construct such a sequence one can begin with $x_n = x$, join $x$ to $e$
by a geodesic segment $\gamma_x$, find the point $x'_{n-1}$ of intersection
of $\gamma_x$ with the sphere $S((n-1)m)$ and a point 
$x_{n-1}\in A_{n-1}\cap B(x'_{n-1}, \epsilon )$, and continue by the 
induction. If $y$ is another point of $A_n$, $(y_0, y_1, \dots , y_n)$
is a corresponding sequence and $k$ is the maximal natural number such
that $x_k = y_k$, then $d(x_{k+1}, y_{k+1})\ge\epsilon$. This
motivates the following definition.

Let us fix $m\in\Bbb N$, $\epsilon > 0$ and $\lambda\in (0,1)$, and denote
by $N_0(n; m, \epsilon , \lambda )$ the maximal cardinality of a subset $A$
of $S(mn)$ satisfying the following condition:

($\ast$) If $x$ and $y$ lie in $A$, then there exist sequences $(x_0, x_1,
\dots x_n)$ and $(y_0, y_1, \dots , y_n)$ of elements of $G$ such that
$x_k,\ y_k\in S(km)$, $x_0 = y_0 = e$, $x_n = x$, $y_n = y$, 
$d(x_j, x_{j+1})\le T = m + \lambda\epsilon$ for all $j$ and 
$d(x_{k+1}, y_{k+1})\ge\epsilon$ when $k$ is the maximal index for which
$x_k = y_k$ (Figure 1).

\midinsert
\centerline{\newbox\rysuubox
\newdimen\rysuuw
\font\rysuua=rysuua at 72.27truept
\setbox\rysuubox=\vbox{\hbox{%
\rysuua\char0\char1\char2\char3\char4\char5}}
\rysuuw=\wd\rysuubox
\setbox\rysuubox=\hbox{\vbox{\hsize=\rysuuw
\parskip=0pt\offinterlineskip\parindent0pt
\hbox{\rysuua\char0\char1\char2\char3\char4\char5}
\hbox{\rysuua\char6\char7\char8\char9\char10\char11}}}
\ifx\parbox\undefined
    \def\setrysuu{\box\rysuubox}
\else
    \def\setrysuu{\parbox{\wd\rysuubox}{\box\rysuubox}}
\fi
\setrysuu
}
\botcaption{Figure 1.}\endcaption
\endinsert

The number
$$\gr^{\rel}(G, S; m, \epsilon , \lambda ) = \limsup_{n\to\infty}\frac{1}{mn}
\log N_0(n; m, \epsilon , \lambda )$$
will be called the exponential {\it rate of growth} of $G$ {\it relative
to} $m$, $\epsilon$ and $\lambda$. Finally, if $\mu :(0,1)\times\Bbb R_+
\to\Bbb N$ and
$\tau :(0,1)\to\Bbb R_+$ are arbitrary functions, then we define the
{\it rate of growth} of $G$ {\it relative to} $\mu$ and $\tau$ by
$$\gr^{\rel}(G, S; \mu , \tau ) =
\sup\{ \gr (G, S; m, \epsilon , \lambda); m > \mu (\lambda ,\epsilon ),
\epsilon > \tau (\lambda ), \lambda\in (0,1)\} .$$

Since $N_0(n; m, \epsilon , \lambda )\le N_0(mn, \epsilon )$ for all
$m$, $n$, $\epsilon$ and $\lambda$, we have
$$\gr^{\rel}(G, S; \mu , \tau )\le\gr (G, S)$$
for all $\mu$ and $\tau$ as above. For the
free group $F_k$ generated by the set $S_k$ of $k$ free generators we have
always
$$\gr^{\rel}(F_k, S_k; \mu, \tau ) = \gr (F_k, S_k).$$
This is because $F_k$ has no "dead ends" (see [{\bf GrH}] for the
definition and some information about some related problems) and
in this case one can arrange $\epsilon$-separated subsets $A_n$ of the
spheres $S(mn)$ in such a way that $\dist (x, A_n) = m$ for any 
$x\in A_{n+1}$. In general, one can expect that a relative
rate of growth is strictly less than the "true" rate of growth.

\subhead 3. Hyperbolic spaces and groups\endsubhead
Let $(X, d)$ be a metric space. A curve $\gamma :[a,b]\to X$ is a {\it
geodesic segment} when 
$$d(\gamma (t), \gamma (s)) = |t - s|$$
for all $t, s\in [a, b]$. The space $X$ is {\it geodesic} when any two
points of $X$ can be joined by a geodesic segment. For any finitely
generated group $G$ and any finite symmetric set $S$ generating $G$, the 
Cayley graph $C(G, S)$ is geodesic.

Given three points $x_0, y$ and $z$ of a metric space $X$, the (based at
$x_0$) {\it Gromov product} of $y$ and $z$ is given by
$$(y|z)_{x_0} = \frac12\left( d(x_0,y) + d(x_0,z) - d(y,z)\right) .$$
The space $X$ is said to be {\it hyperbolic} (more precisely, 
$\delta$-hyperbolic with $\delta\ge 0$) whenever the inequality
$$(x|z)_{x_0}\ge \min\{ (x|y)_{x_0}, (y|z)_{x_0}\} - \delta$$
holds for arbitrary points $x_0, x, y$ and $z$ of $X$. Clearly, the Cayley
graph of any free group $F_k$ ($k = 1, 2, \dots$) generated by the set
$S_k$ of $k$ free generators is a tree, so becomes $0$-hyperbolic. A
finitely generated group is said to be {\it hyperbolic} whenever its
Cayley graph with respect to some (equiv., any) generating set is hyperbolic. 
Free groups and fundamental groups of compact Riemannian manifolds of negative
sectional curvature are hyperbolic.

Assume that $X$ is geodesic, take three points $x_1, x_2$ and $x_3$ of $X$ and
connecting them geodesic segments $\gamma_1, \gamma_2$ and $\gamma_3$. The
union
$$\Delta = \gamma_1\cup\gamma_2\cup\gamma_3$$
is a {\it geodesic triangle} with vertices $x_i$. The triangle $\Delta$
is $\eta$-{\it thin} ($\eta\ge 0$) when the canonical isometry $f_\Delta$ 
mapping $\Delta$ onto a tripod (i.e. the union of three segments with common
origin) $T_\Delta$ satisfies the condition
$$d(x,y)\le d(f_\Delta (x), f_\Delta (y)) + \eta$$
for all $x$ and $y$ of $\Delta$.

In the proof of the Theorem we shall use the following.

\proclaim{Lemma 1} {\rm ([{\bf GhH}], p. 41)} Let $X$ be a geodesic metric
space. If $X$ is $\delta$-hyperbolic, then all the geodesic triangles of
$X$ are $4\delta$-thin. Conversely, if all the geodesic triangles of $X$
are $\eta$-thin, then $X$ is $2\eta$-hyperbolic. \qed
\endproclaim

To construct the boundary $\partial X$ of a hyperbolic space $X$ let us fix a
base point $x_0$ and say that a sequence $(x_n)$ {\it diverges to
infinity} whenever 
$$\lim_{m,n\to\infty}(x_m|x_n) = \infty ,$$
where $(\cdot |\cdot )$ denotes the Gromov product based at $x_0$. Two
such sequences $(x_n)$ and $(y_m)$ are {\it equivalent} whenever
$$\lim_{m,n\to\infty}(x_m|y_n) = \infty .$$
The {\it boundary} $\partial X$ of $X$ consists of all the equivalence classes
of sequences diverging to infinity. Note that $\partial X$ can be
described also in terms of equivalence classes of {\it geodesic rays} 
(i.e., maps $\gamma :[0, \infty)\to X$ such that $\gamma |[0,b]$ is a
geodesic segment for any $b > 0$) or in terms of equivalence classes of
{\it quasirays} (i.e. quasi-isometric maps of $[0,\infty )$ into $X$): Two
such rays (or, quasirays) $\gamma$ and $\sigma$ are equivalent whenever their
Hausdorff distance $d_H(\gamma , \sigma )$ is finite. The boundary point
corresponding to the equivalence class of such $\gamma$ is that determined
by the sequence $x_n = \gamma (n)$, $n\in\Bbb N$. The equivalence of these
constructions follows from the following fact which will be used later.

\proclaim{Lemma 2} {\rm ([{\bf GhH}], p. 87)} Let $X$ be a
$\delta$-hyperbolic geodesic metric space. For any $c \ge 1$ there exists
$D\ge 0$ such that any $c$-quasigeodesic segment $\gamma$ (i.e. any
map $\gamma :[a, b]\to X$ (resp., $\gamma :[a,b]\cap\Bbb Z\to X$) such that 
the inequality $- c - c|s - t|\le d(\gamma (s), \gamma (t))\le c|s - t| + c$ 
holds for all $s$ and $t$ of $[a, b]$ (resp., of $[a,b]\cap\Bbb Z$)) and any 
geodesic segment $\sigma$ joining $\gamma (a)$ to $\gamma (b)$ satisfy the 
inequality
$$d_H(\gamma , \sigma )\le D. \qed$$
\endproclaim
The set $\partial X$ can be equipped with the metric structure as follows.
First, for any $\xi$ and $\zeta$ of $\partial X$ put
$$(\xi |\zeta ) = \sup\liminf_{m, n\to\infty}(x_m|y_n),$$
where $(x_m)$ and $(y_n)$ run over the set of all diverging to infinity
sequences representing, respectively, $\xi$ and $\zeta$.
Note that if $X$ is $\delta$-hyperbolic, then
$$(\xi |\zeta )- 2\delta\le\liminf_{m, n\to\infty}(x_m|y_n)\le (\xi |\zeta
)$$
for all sequences $(x_m)$ and $(y_m)$ representing $\xi$ and $\zeta$. 
Next, choose $\eta > 0$ and put
$$\rho_\eta (\xi , \zeta ) = \exp (-\eta\cdot (\xi |\zeta )).$$
Finally, let
$$d_\eta (\xi , \zeta ) = \inf\{\sum_{i=0}^k\rho_\eta (\xi_i, \xi_{i+1});
\xi_i\in\partial X, \xi_0 = \xi \ \text{and}\ \xi_{k+1} = \zeta , k\in\Bbb N\}.
$$
If $\eta > 0$ is small enough, then $d_\eta$ is a distance function
on $\partial X$ and $(\partial X, d_\eta )$ becomes a compact metric space
of finite Hausdorff dimension ([{\bf GhH}], pp. 122 - 126). Moreover, the 
inequalities 
$$(1 - 2\eta ')\rho_\eta (\xi , \zeta )\le d_\eta (\xi , \zeta )\le
\rho_\eta (\xi , \zeta ), \quad \xi , \zeta\in\partial X,$$
hold with 
$$\eta ' = \exp (\eta\delta ) - 1.$$

\subhead 4. Proof of the Theorem\endsubhead
Let again $G$ be a finitely generated group, $S$ - a finite symmetric set
generating $G$ and consider $C(G, S)$, the Cayley graph of $G$ equipped with the
distance function $d$ satisfying
$$d(g_1, g_2) = |g_1^{-1}g_2|,\quad g_1, g_2\in G,$$
where $|\cdot |$ is the length function on $G$ determined by $S$, and
making the edges of $C(G,S)$ isometric to Euclidean segments of length 1. 
Assume that $C(G, S)$ is $\delta$-hyperbolic and let $\partial G = 
\partial C(G, S)$ be its boundary equipped, as in Section 3, with the distance
function $d_\eta$, $\eta > 0$ being small enough. 

The group $G$ acts on $C(G, S)$ {\it via} isometries which, when restricted to
$G\subset C(G, S)$, reduce to left translations $L_g$, $G\ni h\mapsto gh$.
Therefore, each $L_g$ extends to a Lipschitz homeomorphism, denoted by
$L_g$ again, of the boundary $\partial G$. Given $\theta > 0$ denote by $N(n, \theta ;
\partial G)$ the maximal number of points of $\partial G$ pairwise 
$(n, \theta )$-separated by this action of $G$.  Also, let $N'(n,
\theta ; \partial G)$ be the minimal cardinality of an $(n, 
\theta )$-spanning subset of $\partial G$ and $h(G, S; \partial G)$ - the
corresponding entropy. 

Since $G$ is hyperbolic, 
there exists a constant $c_0 > 1$ such that for any $g\in G$ there exists 
$g'\in G$ such that $d(g,g')\le c_0$ and $|g'| = |g| + 1$ (see [{\bf GrH}],
p. 60). Let $D$ be a corresponding constant such that any $c_0$-quasigeodesic
segment  lies at most $D$-apart (in the Hausdorff distance) from a true 
geodesic (compare Lemma 2). Let
$$\tau (\lambda ) = \frac{4(D+\delta )}{1 - \lambda}\quad\text{and}\quad
\mu (\lambda , \epsilon ) = \frac{\lambda\epsilon}{c_0 - 1}.$$
We shall show that

$$(\ast\ast )\qquad \gr^{\rel}(G, S; \mu , \tau )
\le h(G, S; \partial G)\le\gr (G, S).$$

We begin by the proof of the second inequality in $(\ast\ast )$.

Choose $\theta > 0$ and a natural number $k$ for which the inequality 
$$\exp (-\eta\cdot k) < \theta$$
holds. For any $n\in\Bbb N$ and any $g\in S(n+k)$ choose, if only possible, 
a point $\xi_g\i\partial G$ which can be connected to $e$ by a geodesic ray, 
say $\gamma_g$, which passes through $g$. We claim that the set
$$A_n = \{\xi_g ; g\in S(n+k)\}$$ 
is $(n, \epsilon )$-spanning in $\partial G$. Indeed, if $\zeta\in\partial G$,
$\sigma :[0, \infty )\to X$ is a geodesic ray connecting $e$ to $\zeta$, $h\in
S(n)$ and $g = \sigma (n+k)$, then $g\in S(n+k)$, the corresponding ray
$\gamma_g$ exists and satisfies the conditions
$$d(x_j, y_i)\le i + j - 2(n+k)$$
and 
$$2(h^{-1}y_i|h^{-1}x_j)\ge (i-n) + (j-n) - (i+j) + 2(n+k) = 2k,$$
where $x_j = \sigma (j)$ and $y_i = \gamma (i)$ for all $i$ and $j$ 
sufficiently large (Figure 2). 

\midinsert
\centerline{\newbox\rysubox
\newdimen\rysuw
\font\rysua=rysua at 72.27truept
\setbox\rysubox=\vbox{\hbox{%
\rysua\char0\char1\char2\char3\char4\char5}}
\rysuw=\wd\rysubox
\setbox\rysubox=\hbox{\vbox{\hsize=\rysuw
\parskip=0pt\offinterlineskip\parindent0pt
\hbox{\rysua\char0\char1\char2\char3\char4\char5}
\hbox{\rysua\char6\char7\char8\char9\char10\char11}}}
\ifx\parbox\undefined
    \def\setrysu{\box\rysubox}
\else
    \def\setrysu{\parbox{\wd\rysubox}{\box\rysubox}}
\fi
\setrysu
}
\botcaption{Figure 2.}\endcaption
\endinsert

Therefore,
$$d_\eta (L_{h^{-1}}\zeta , L_{h^{-1}}\xi_g) \le
\rho_\eta (L_{h^{-1}}\zeta , L_{h^{-1}}\xi_g) < e^{-k\eta} < \theta .$$
This shows the inequalities
$$N'(n, \theta ; \partial G)\le\# A_n\le\# S(n+k)\le N(n+k)$$
which imply immediately the required inequality in $(\ast\ast )$.

The proof of the first inequality in $(\ast\ast )$ is a bit more complicated.

Fix $\lambda\in (0,1)$, $\epsilon > \tau (\lambda )$ and
$m > \mu (\lambda , \epsilon )$. For any $n\in\Bbb N$ choose a maximal
subset $A_n$ of $S(mn)$ satisfying condition ($\ast$) of Section 2. 
For any $x\in A_n$ set $x_n = x$, choose a point 
$x_{n-1}\in A_{n-1}$ such that $d(x, x_{n-1})\le m + \lambda\epsilon$, then 
a point $x_{n-2}\in A_{n-2}$ for which $d(x_{n-1}, x_{n-2})\le m + 
\lambda\epsilon$ and so on. Finally, put $x_0 = e$. The map
$$\{ 0, m , \dots , mn\}\ni j\mapsto x_{j/m}$$
is  $c$-quasi-isometric with $c = (m + \lambda\epsilon )/m < c_0$.

Each map considered above can be extended to a $c_0$-quasi-isometric map 
$$\gamma_x: \Bbb N\ni j\mapsto x'_j\in G$$ such that $x'_{im} = x_i$ for
$i = 0,1,\dots , n$ and the sequence $(x_j)$ converges to a point $\xi_x$ of
$\partial G$. We are going to show that  the set
$$\{ \xi_x; x\in A_n\}$$
is $(n, \theta )$-separated under the action of $G$ for some $\theta$
independent of $n$.

To this end, let us take arbitrary points $x$ and $y$ of $A_n$, $x\ne y$,
choose sequences $(x_0, x_1, \dots , x_n)$ and $(y_0, y_1, \dots , y_n)$
as above, and let $k$ be the maximal element of $\{ 0, 1, \dots , n\}$ for which
$x_k = y_k$. Denote this common value of $x_k$ and $y_k$ by $g$ and consider the
quasi-geodesic rays $\bar\gamma_x$ and $\bar\gamma_y$ obtained from $\gamma_x$
and $\gamma_y$ by restricting their domains to $\{ mk, mk+1, \dots\}$. Then, 
$\bar\gamma_x$ and $\bar\gamma_y$ originate at $g$ and converge to $\xi_x$
and $\xi_y$, respectively. By Lemma 2, there exist geodesic rays 
$\tilde\gamma_x$ and $\tilde\gamma_y$ originated at $e$ and within 
the Hausdorff distance $D$ from $L_{g^{-1}}\circ\bar\gamma_x$ and
$L_{g^{-1}}\circ\bar\gamma_y$, respectively. Then, for
any $j\in\Bbb N$ there exist positive real numbers $s_j$ and $t_j$ such that
$$d(\tilde\gamma_x (s_j), g^{-1}\bar\gamma_x (j))\le D\quad\text{and}\quad
d(\tilde\gamma_y (t_j), g^{-1}\bar\gamma_y (j))\le D.$$
In particular,
$$d(\tilde\gamma_x (s_m), g^{-1}x_{k+1})\le D\quad\text{and}\quad
d(\tilde\gamma_y (t_m), g^{-1}y_{k+1})\le D.$$
Clearly,  
$$m - D\le s_m, t_m\le m + D + \lambda\epsilon ,\  |s_m - t_m|\le 2D +
\lambda\epsilon$$
and
$$d(\tilde\gamma_x (s_m), \tilde\gamma_y (t_m))\ge\epsilon - 2D.$$
Assume that $(\tilde\gamma_x (s_j)|\tilde\gamma_y(t_j)) > m + D +
\lambda\epsilon$ for some $j\in\Bbb N$. The isometry $f_\Delta$ corresponding
to the geodesic triangle $\Delta$ with vertices $e, \tilde\gamma_x (s_j),
\tilde\gamma_y (t_j)$ (compare Lemma 1) maps the points $\tilde\gamma_x (s_m)$
and $\tilde\gamma_y (t_m)$ onto some points of the originated at $e$ edge
of the tripod $f_\Delta (\Delta )$  and, therefore, satisfies the condition
$$\align d(\tilde\gamma_x(s_m), \tilde\gamma_y(t_m))&\le 
d(f_\Delta (\tilde\gamma_x(s_m)), f_\Delta (\tilde\gamma_y(t_m)))
+ 4\delta\\ &= |s_m - t_m| + 4\delta \le 2D + \lambda\epsilon + 4\delta
< \epsilon - 2D.
\endalign$$
Comparing the inequalities above we obtain a contradiction which shows
that
$$(\tilde\gamma_x (s_j)| \tilde\gamma_y(t_j))\le m + D + \epsilon\lambda$$
for all $j\in\Bbb N$. This inequality proves that
$$\align
d_\eta (g^{-1}\xi_x, g^{-1}\xi_y)&\ge 
(1 - 2\eta ')\rho_\eta (g^{-1}\xi_x, g^{-1}\xi_y)\\ &\ge 
(1- 2\eta ')\exp (-\eta (m + D + \lambda\epsilon + 2\delta )),
\endalign$$
i.e. that the points $\xi_x$ and $\xi_y$ are $(nm, \theta )$-separated with
$$\theta = (1 - 2\eta ')\exp (-\eta (m + D + \lambda\epsilon + 2\delta)).$$

The above argument implies the inequality
$$N(nm, \theta ; \partial G)\ge\# A_n\ge N_0(G,S; m, \epsilon , \lambda ).$$
which holds for all $n$. Passing to suitable limits when $n\to\infty$ yields
the required inequality in ($\ast\ast$). \qed

In [{\bf Fr}], Friedland defined the {\it minimal entropy} $h_{\min}(G)$
of a finitely generated group $G$ of homeomorphisms of a compact metric space
$X$:
$$h_{\min}(G) = \inf_{S}h(G, S),$$
where $S$ ranges over all finite symmetric sets generating $G$. Similarly,
the {\it minimal rate of growth} $\gr_{\min}(G)$ of any finitely generated 
group $G$ can be defined as follows (compare [{\bf GrH}]):
$$\gr_{\min}(G) = \inf_S\gr (G, S).$$
The reader can define the minimal relative rate of growth
$\gr_{\min}^{\rel}(G)$ appropriately.

If $G$ is hyperbolic, then $\id_G$ induces a H\"older homeomorphism of 
boundaries of $G$ obtained from different generating sets ([{\bf GhH}],
page 128). Therefore, the boundary entropy of such $G$ (w.r.t. a given
finite symmetric generating set $S$) does not depend on the choice of
a generating set used in the construction of $\partial G$ and our Theorem
implies immediately the following.

\proclaim{Corollary} For any hyperbolic group $G$ the equalities
$$\gr_{\min}^{\rel}(G)\le h_{\min}(G, \partial G)\le\gr_{\min}(G)$$
hold. \qed
\endproclaim

This answers partially the following question asked by Friedland in 
[{\bf Fr}]: Find a geometric interpretation of the minimal entropy of 
a Kleinian group acting on the ideal boundary of a hyperbolic space H$^n$.

For the free group $F_k$ with $k$ generators the above discussion and an
argument of [{\bf GLP}] (p. 70) imply the equality
$$h_{\min}(F_k, \partial F_k) = \gr_{\min}(F_k) = \gr_{\min}^{\rel} (F_k) =
\log (2k-1).$$
In fact, if $S$ is any finite symmetric set generating $F_k$, then the 
elements of $S$ represent members of a set $S'$ generating $\Bbb Z^k$, the 
abelianization of $F_k$. $S'$ contains a symmetric set $R'$ such that 
$\# R' = 2k$ and the subgroup of $\Bbb Z^k$ generated by $R'$ has finite 
index. The corresponding subset $R$ of $S$ consists also of $2k$ elements and
generates the free group isomorphic to $F_k$. Therefore,
$$\align h(F_k, S; \partial F_k)&\ge h(F_k, R; \partial F_k) = h(F_k,
S_k; \partial F_k)\\ & \ge\gr^{\rel} (F_k, S^k) = \gr (F_k, S_k) = \log
(2k-1).\endalign$$
The opposite inequality is obvious.

The other natural case, that of the fundamental group $\Gamma_g$ of a closed 
oriented surface of genus $g > 1$ is more complicated: The minimal rate of 
growth of $\Gamma_g$ is still unknown even if some estimates exist: 
$\gr_{\min}(\Gamma_g)\ge 4g - 3$, $\gr_{\min}(\Gamma_g)\le\gr (\Gamma_g, S_g)
\approx 4g - 1 - \epsilon_g$, where $S_g$ is the canonical set of generators
of $\Gamma_g$ and $\epsilon_g$ is a small constant found numerically.
In particular, $5\le\gr_{\min}(\Gamma_2)\le 6.9798$. The calculation or
estimation of the value of the minimal relative rate of growth of $\Gamma_g$
is yet more difficult.


\Refs 
\widestnumber\key{\bf GLW} 

\ref\key{\bf CP}\by M. Coornaert and A. Papadopoulos\book Symbolic
Dynamics and Hyperbolic Groups\publ Springer Verlag\publaddr
Berlin - Heidelberg\yr 1993\endref

\ref\key{\bf Eg}\by S. Egashira\paper Expansion growth of foliations\jour Ann. Fac.  
Sci. Univ. Toulouse\vol 2\yr 1993\pages 15--52\endref 

\ref\key{\bf Fr}\by S. Friedland\paper Entropy of graphs, semigroups and 
groups\inbook Ergodic theory of $\Bbb Z^d$ actions\eds M. Policott and 
K. Schmidt\publ London Math. Soc.\bookinfo {Lecture Notes Ser., vol. {\bf 228}}
\publaddr London\yr 1996\pages 319 -- 343\endref

\ref\key{\bf GhH}\by E. Ghys and P. de la Harpe\book Sur les Groupes 
Hyperboliques d'apres Mikhael Gromov\publ Birkh\"auser\publaddr
Boston - Basel - Berlin \yr 1990\endref

\ref\key{\bf GLW}\by E. Ghys, R. Langevin and P. Walczak\paper Entropie  
g\'eom\'etrique des feuilletages\jour Acta Math.\vol 160\yr 1988
\pages 105--142\endref 

\ref\key{\bf GrH}\by R. Grigorchuk and P. de la Harpe\paper On problems
related to growth, entropy and spectrum in group theory\jour J. Dyn.
Control Sys.\vol 3\yr 1997\pages 51--89\endref

\ref\key{\bf Gro}\by M. Gromov\paper Hyperbolic groups\inbook Essays in
group theory\eds S. M. Gersten\publ Springer
Verlag\publaddr Berlin -- Heidelberg -- New York\yr 1987
\pages 75 -- 263\endref

\ref\key{\bf GLP}\by M. Gromov, J. Lafontaine and P. Pansu\book
Structures m\'etriques pur les vari\'et\'es riemanniennes\publ Cedic/F.
Nathan\publaddr Paris\yr 1981\endref

\ref\key{\bf Wa}\by P. Walters\book An introduction to Ergodic Theory\publ
Springer Verlag\publaddr New York - Heidelberg - Berlin\yr 1982\endref
\endRefs
\enddocument